\documentclass[10pt]{amsart}
\usepackage{amsthm,amsmath,amssymb,amsfonts}
\usepackage[american]{babel}
\usepackage{cite}
\usepackage{url}

\allowdisplaybreaks

\newtheorem{theorem}{Theorem}[section] 
\newtheorem*{theorem*}{Theorem}

\newtheorem{proposition}[theorem]{Proposition}

\theoremstyle{definition}

\theoremstyle{remark}
\newtheorem{remark}[theorem]{Remark}

\numberwithin{equation}{section}

\begin{document}


\newcommand{\abs}[1]{\lvert#1\rvert}
\newcommand{\NU}[1]{\mbox{\rm V}(#1)}
\newcommand{\NUreg}[1]{\mbox{\rm V}(#1)_{\text{\rm reg}}^{\text{\rm tor}}}
\newcommand{\UNIT}[1]{\mbox{\rm U}(#1)}
\newcommand{\cen}[2]{\mbox{\rm C}_{#1}(#2)}
\newcommand{\nor}[2]{\mbox{\rm N}_{#1}(#2)}
\newcommand{\zen}[1]{\mbox{\rm Z}(#1)}
\newcommand{\lara}[1]{\langle{#1}\rangle}
\newcommand{\ZZ}{\mbox{$\mathbb{Z}$}}
\newcommand{\sZZ}{\mbox{$\scriptstyle\mathbb{Z}$}}
\newcommand{\QQ}{\mbox{$\mathbb{Q}$}}
\newcommand{\sQQ}{\mbox{$\scriptstyle\mathbb{Q}$}}
\newcommand{\FF}{\mbox{$\mathbb{F}$}}
\newcommand{\sFF}{\mbox{$\scriptstyle\mathbb{F}$}}
\newcommand{\NN}{\mbox{$\mathbb{N}$}}
\newcommand{\CC}{\mbox{$\mathbb{C}$}}
\newcommand{\paug}[2]{\varepsilon_{#1}(#2)}


\title[Partial augmentations and Brauer character values]
{Partial augmentations and Brauer character values
of torsion units in group rings}

\author{Martin Hertweck}
\address{Universit\"at Stuttgart, Fachbereich Mathematik,
Institut f\"ur Geometrie und Topo\-logie,    
Pfaffenwald\-ring 57, 70550 Stuttgart, Germany}
\email{hertweck@mathematik.uni-stuttgart.de}

\thanks{The author acknowledges 
support by the Deutsche Forschungsgemeinschaft.}
\date{\today. \\
Revised version of manuscript from December 31, 2004 
({\tt arXiv:math.RA/0612429})}


\begin{abstract}
For a torsion unit $u$ of the integral group ring $\ZZ G$ of a finite 
group $G$, and a prime $p$ which does not divide the order of $u$ 
(but the order of $G$), a relation between the partial augmentations 
of $u$ on the $p$-regular classes of $G$ and Brauer character values 
is noted, analogous to the obvious relation between partial augmentations
and ordinary character values.
For non-solvable $G$, consequences concerning rational conjugacy of $u$ 
to a group element are discussed, considering as examples the
symmetric group $S_{5}$ and the groups $\text{\rm PSL}(2,p^{f})$.
\end{abstract}

\maketitle

\section{Introduction}\label{Sec:intro}

This note is about an observation on torsion units in the integral group ring 
$\ZZ G$ of a finite group $G$ which seemingly has been overlooked so far. 
Sure, we can restrict our attention to units of augmentation (= sum 
of coefficients) one of $\ZZ G$, which form a group we denote by $\NU{\ZZ G}$.
A torsion unit in $\NU{\ZZ G}$ is said to be rationally 
conjugate to a group element if it is conjugate to an element of $G$
by a unit of the rational group ring
$\QQ G$, and according to a conjecture of Hans Zassenhaus, 
this should always be the case. For the present state of affairs, we refer 
the reader to \cite{Seh:93,Seh:02,He:04a,He:06a}.

A close connection with partial augmentations was noted in 
\cite{MaRiSeWe:87}: a torsion unit  $u$ in $\NU{\ZZ G}$ is rationally 
conjugate to a group element if and only if for each power of $u$,
all but one of its partial augmentations vanish. 
Recall that for a group ring element $u=\sum_{g\in G}a_{g}g$ 
(all $a_{g}$ in $\ZZ$), its partial augmentation
with respect to an element $x$ of $G$, or rather its
conjugacy class $x^{G}$ in $G$, is the sum $\sum_{g\in x^{G}}a_{g}$;
we will denote it by $\paug{x}{u}$.
Further connections are described in Section~\ref{Sec:PA}, 
but let us mention another known one
which is of vital importance in this work.
Let $R$ be a Dedekind ring of characteristic zero in which a rational 
prime $p$ is not invertible, and let $u$ be a $p$-regular torsion unit 
in $\NU{RG}$ (the given definitions also make sense for more general 
coefficient rings). Then $\paug{x}{u}=0$ for every $p$-singular element
$x$ of $G$. (A torsion unit is called $p$-regular if $p$
does not divide its order, and $p$-singular otherwise.) 

Brauer's definition of $p$-modular characters (now called Brauer characters)
as certain complex-valued functions $\varphi$ on $G_{\text{\rm reg}}$ 
($=$ set of $p$-regular elements of $G$) can naturally be extended to make 
$\varphi(u)$ meaningful for elements $u$ of the set $\NUreg{RG}$ of 
$p$-regular torsion units in $\NU{RG}$, and then a few elementary facts
about Brauer characters still hold (see Proposition~\ref{P1}).
From that, we deduce our main observation (Theorem~\ref{MT}):
\[ \varphi(u)=
\sum_{\substack{x^{G}:\; x\text{ is}\\ \text{$p$-regular}}}
\varepsilon_{x}(u)\varphi(x)
\quad \text{for all $u\in\NUreg{RG}$.} \]

This allows us to formulate (in Section~\ref{Sec:LPmethod})
a $p$-modular version of a method introduced by  Luthar and Passi 
\cite{LuPa:89} which imposes constraints on the partial augmentations 
of $p$-regular torsion units in $\ZZ G$.
This method may prove useful in investigating the Zassenhaus conjecture
for non-solvable groups. 
In Section~\ref{Sec:ex1}, it is used to verify the Zassenhaus conjecture 
for the symmetric group $S_{5}$; based on a computer calculation, this
was done before by Luthar and Trama \cite{LuTr:91}.
In Section~\ref{Sec:ex2}, it is used to explore the Zassenhaus conjecture 
for the groups $\text{\rm PSL}(2,p^{f})$.
Some evidence is given that the conjecture might by verified for the 
groups $\text{\rm PSL}(2,p)$; a complete proof is given for $p=7,11,13$.

We remark that the Luthar--Passi method is of computational nature,
and it seems worthwhile to examine more examples on a computer 
to see how far it will take us.
The computer algebra system GAP \cite{GAP4} lends itself to this task.
It makes available the tables contained in the Atlas of Brauer Characters
\cite{JLPW:95} and many others, and provides functions for computations
with characters and character tables.
Since the first version of this paper was written, 
this line of research has been pursued in
\cite{BoKo:06,BoKo:06b,BoJeKo:06,BoKoSi:06,BoHe:06}.

\section{Partial augmentations}\label{Sec:PA}

We collect some facts about partial augmentations of torsion units in 
integral group rings which will be used later on.
Only Proposition~\ref{lemmi} is mathematically new in this section,
but has a precursor in \cite[Proposition~3.1]{He:04a}.
Fix a finite group $G$, and let $R$ be an integral domain of 
characteristic zero, with quotient field $K$. 

Our tool to prove rational conjugacy of torsion units is 
\cite[Theorem~2.5]{MaRiSeWe:87} (formulated in \cite[Lemma~2.5]{He:04a} 
for more general coefficient rings) which, because of its fundamental
importance, shall be reproduced here.
\begin{theorem}\label{paMRSW}
Let $u$ be a torsion unit in $\NU{RG}$, with no prime divisor of the
order of $u$ being invertible in $R$. Then $u$ is conjugate to an element 
of $G$ by a unit of $KG$ if and only if for every $m$ dividing the order
of $u$, all partial augmentations of $u^{m}$ but one vanish. 
\end{theorem}

Vanishing of partial augmentations can be caused for arithmetical reasons.
A classical result, due to Berman and Higman (see \cite[(1.4)]{Seh:93}), 
states that a torsion unit in $\NU{\ZZ G}$ is a central element of $G$
provided that its partial augmentation with respect to such an element
is nonzero.

It is not known, even for solvable $G$, whether for a given torsion unit $u$
in $\NU{\ZZ G}$, there exists an element in $G$ of the same order, possibly
such that the associated partial augmentation is nonzero
(Research Problem~8 from \cite{Seh:93}).
The answer is affirmative for units $u$ of prime power order,
by \cite[Theorem~4.1]{CoLi:65}. (The used technique
is also described in \cite[\S7]{Seh:93}, where the result is 
attributed to Zassenhaus.) Moreover, if the order of $u$ is $r^{n}$,
for a prime $r$, then
\begin{align}\label{CLequi}
\sum_{\substack{x^{G}:\; x\text{ of}\\ \text{order $r^{m}$}}}
\varepsilon_{x}(u) \equiv 0 \mod{r} \quad\text{for all $m<n$.}
\end{align}
That these congruences actually might be equalities is known as a 
conjecture of A.~A.~Bovdi (rendered in \cite{Jur:95}). For
$G=\text{\rm PSL}(2,p^{f})$, it is confirmed in Proposition~\ref{Bov}.

The following proposition is an obvious generalization of 
\cite[Proposition~3.1]{He:04a}, and can be proved in the exact same manner.
Previously, only special cases were known
(i.e., \cite[Theorem~2.7]{MaRiSeWe:87}).
\begin{proposition}\label{lemmi}
Let $R$ be a Dedekind ring of characteristic zero in which a given
rational prime $p$ is not invertible, and let $u$ be a torsion unit 
in $\NU{RG}$. Then $\paug{x}{u}=0$ for every element $x$ of $G$ whose 
$p$-part has strictly larger order than the $p$-part of $u$. 
\end{proposition}
\begin{proof}
Suppose that $x$ is an element of $G$ whose $p$-part has strictly larger
order than the $p$-part of $u$. 
Let $C$ be an (abstract) cyclic group whose order is the least common 
multiple of the orders of $u$ and $x$.
Enlarging $R$, if necessary, we can assume that $R$ is a complete 
discrete valuation ring (cf.\ \cite[\S4c]{CuRe:81}).
Set $M=RG$, viewed as $RC$-lattice by letting a generator $c$ 
of $C$ act by $m\cdot c=x^{-1}mu$ for $m\in M$.
By \cite[(38.12)]{Seh:93}, we have to show that $\chi(c)=0$ for the 
character $\chi$ of $C$ afforded by $KM$.
Note that $c$ is $p$-singular, so this will follow once we have shown
the stronger statement that $M$ is a projective $RC$-lattice, by
Green's Theorem on Zeros of Characters (see \cite[(19.27)]{CuRe:81}).
This in turn follows from the assumption, meaning that 
the action of the subgroup $P$ of order $p$ in $C$ is given
by a multiplication action of the subgroup of order $p$ in $\lara{x}$,
which shows that $M$ is a free $RP$-lattice.
We provide details.
Let $\pi$ be a prime element of $R$, and set $k=R/\pi R$ 
(a residue field of characteristic $p$). It is enough 
to show that $k\otimes_{R}M$ as $kC$-module is projective
(see \cite[(20.10), (30.11)]{CuRe:81}), i.e., that
$k\otimes_{R}M$ is projective relative to a Sylow $p$-subgroup of $C$.
It is well-known that this follows from the projectivity of
$k\otimes_{R}M$ as $kP$-module (see \cite[Lemma~3.2]{He:04a}).
\end{proof}

This applies in particular to torsion units in $\ZZ G$, giving an 
affirmative answer to Research Problem~9 from \cite{Seh:93}:
\begin{theorem}\label{T8}
Let $u$ be a torsion unit in $\NU{\ZZ G}$. Then $\paug{x}{u}\neq 0$ is 
possible only for elements $x$ of $G$ whose order is a divisor of the 
order of $u$. \qed
\end{theorem}

\section{Brauer character values}\label{Sec:MC}

A general reference for the following discussion is Chapter~2 
(in particular \S17) of \cite{CuRe:81}.
Let $G$ be a finite group, $p$ a prime divisor of its order,
and fix a $p$-modular system $(K,R,k)$ such that
$\text{\rm char}(K)=0$, with $K$ sufficiently large relative to $G$.
Then also $k=R/\pi R$ is sufficiently large relative to $G$.
(Here, $\pi$ is a prime element of the discrete valuation ring $R$.)

Brauer associated with a modular representation 
$\Theta:G\rightarrow\text{\rm GL}(n,k)$ a complex-valued function $\varphi$
on $G_{\text{\rm reg}}$ ($=$ set of $p$-regular elements of $G$), now
called a Brauer character of $G$.
The way Brauer did this actually shows that 
we can extend the domain of $\varphi$ to the set $\NUreg{RG}$ of 
$p$-regular torsion units in $\NU{RG}$.
Having the natural map $R\rightarrow R/\pi R=k$ at our disposal, 
we consider the representation $\Theta$ as an algebra homomorphism
from $RG$ into the algebra of $n\times n$-matrices over $k$. 
To define $\varphi(u)$ for a unit $u$ in $\NUreg{RG}$,
notice that the eigenvalues of $\Theta(u)$ are certain roots of unity
in $k$, of order relatively prime to $p$, and their sum is the trace of
$\Theta(u)$ (the usual character). We just replace these eigenvalues by
the roots of unity in $R$ to which they are in bijection
via the map $R\rightarrow k$, and call the resulting sum $\varphi(u)$.

We fix elements $x_{1},\dotsc,x_{r}$ of $G$
which are representatives of the $p$-regular conjugacy classes of $G$.
If $u$ is a $p$-regular torsion unit in $\NU{RG}$, then 
$\varepsilon_{y}(u)=0$ for every $p$-singular element $y$ of $G$
(see Proposition~\ref{lemmi}).
Hence for an ordinary $K$-character $\chi$ of $G$,
\begin{align}\label{Eq2}
\chi(u)=\sum_{l=1}^{r}\varepsilon_{x_{l}}(u)\chi(x_{l})\quad
\text{for all $u\in\NUreg{RG}$.}
\end{align}
We shall see in a moment that the analogue holds for Brauer characters.

The following elementary facts about these (extended) Brauer characters are 
proved in the exact same manner as \cite[(17.5)]{CuRe:81}.
\begin{proposition}\label{P1}
Let $u\in\NUreg{RG}$. Then the following holds.
\begin{enumerate}
\item[(i)]
Let $\varphi$ be the Brauer character of a $kG$-module $M$.
Let $L$ be a submodule of $M$, let $\phi_{L}$ be the Brauer character 
afforded by $L$, and let $\phi_{M/L}$ be the Brauer character 
afforded by $M/L$. Then $\varphi(u)=\phi_{L}(u)+\phi_{M/L}(u)$.
\item[(ii)]
Let $V$ be a $KG$-module with $K$-character $\chi$, and let $\varphi$
be the Brauer character of $M/\pi M$, where $M$ is a full $RG$-lattice in $V$.
Then $\chi(u)=\varphi(u)$.
\end{enumerate}
\end{proposition}

Let $\chi_{1},\dotsc,\chi_{h}$ be the ordinary irreducible $K$-characters of $G$,
and let $\varphi_{1},\dotsc,\varphi_{r}$ be the irreducible 
Brauer characters of $G$.
Let $\mathbf{D}=(d_{ij})$ be the decomposition matrix of $G$
(relative to $p$).
By its definition, and Proposition~\ref{P1}, we have
\begin{align}\label{Eq1}
\chi_{i}(u)=\sum_{j=1}^{r}d_{ij}\varphi_{j}(u)\quad
\text{for all $u\in\NUreg{RG}$.}
\end{align}

We now prove the main result of this note.
\begin{theorem}\label{MT}
Let $\varphi$ be a Brauer character of $G$ (relative to $p$), and
let $u$ be a $p$-regular torsion unit in $\NU{RG}$. Then,
with representatives $x_{1},\dotsc,x_{r}$ of the $p$-regular 
conjugacy classes of $G$,
\[ \varphi(u)=\sum_{l=1}^{r}\varepsilon_{x_{l}}(u)\varphi(x_{l}). \]
\end{theorem}
\begin{proof}
By Proposition~\ref{P1}(i), it is enough to prove the result for an
irreducible Brauer character $\varphi_{j}$.
By (\ref{Eq1}) and (\ref{Eq2}), we have
\begin{align*}
\mathbf{D}
\left(\begin{matrix}
\varphi_{1}(u) \\ \vdots \\ \varphi_{r}(u) 
\end{matrix} \right) & =
\left(\begin{matrix}
\chi_{1}(u) \\ \vdots \\ \chi_{h}(u) 
\end{matrix} \right) =
\sum_{l=1}^{r}\varepsilon_{x_{l}}(u)
\left(\begin{matrix}
\chi_{1}(x_{l}) \\ \vdots \\ \chi_{h}(x_{l}) 
\end{matrix} \right) =
\sum_{l=1}^{r}\varepsilon_{x_{l}}(u)\;
\mathbf{D}
\left(\begin{matrix}
\varphi_{1}(x_{l}) \\ \vdots \\ \varphi_{r}(x_{l}) 
\end{matrix} \right) 
\\
& = \mathbf{D}
\left(\begin{matrix}
\sum_{l=1}^{r}\varepsilon_{x_{l}}(u)\varphi_{1}(x_{l}) \\ \vdots \\
\sum_{l=1}^{r}\varepsilon_{x_{l}}(u)\varphi_{r}(x_{l}) 
\end{matrix} \right).
\end{align*}
The $r\times r$-matrix $\mathbf{D}^{\text{\rm t}}\mathbf{D}$ is the 
Cartan matrix $\mathbf{C}$, and $\det(\mathbf{C})\neq 0$. Hence
\[ \left(\begin{matrix}
\varphi_{1}(u) \\ \vdots \\ \varphi_{r}(u) 
\end{matrix} \right) =
\left(\begin{matrix}
\sum_{l=1}^{r}\varepsilon_{x_{l}}(u)\varphi_{1}(x_{l}) \\ \vdots \\
\sum_{l=1}^{r}\varepsilon_{x_{l}}(u)\varphi_{r}(x_{l}) 
\end{matrix} \right), \]
and we are done.
\end{proof}

\begin{remark}\label{remi2}
Suppose that $u$ is a $p$-regular torsion unit in $\NU{RG}$ which is 
conjugate to an element $g$ of $G$ by a unit of $KG$. Then $u$ is even 
conjugate to $g$ by a unit of $RG$ (see \cite[Lemma~2.9]{He:04a}), whence
the images of $u$ and $g$ in $kG$ are conjugate by a unit of $kG$. 
In particular, if $\Theta$ is a $p$-modular representation of $G$,
then $\Theta(u)$ has the same eigenvalues as $\Theta(g)$.
\end{remark}

\section{The Luthar--Passi method}\label{Sec:LPmethod}

The method introduced by  Luthar and Passi \cite{LuPa:89} imposes
constraints---once ordinary characters of the finite group $G$ are 
known---on the partial augmentations of a torsion unit in $\ZZ G$.
The idea behind is to exploit the simple fact that the 
character values of a torsion unit, which are determined by its partial 
augmentations and the character table, are sums of roots of unity.
It is quite obvious that Theorem~\ref{MT} allows us to formulate, for 
$p$-regular torsion units in $\ZZ G$, a $p$-modular version of this method,
which can be applied once Brauer characters of $G$ are known.

Let $(K,R,k)$ be a $p$-modular system as in Section~\ref{Sec:MC}, 
and let the field $F$ be either $K$ or $k$.
Let $U$ be an invertible matrix with entries in $F$ such that
$U^{n}$ is the identity matrix for some natural number $n$.
(In a moment, $U$ will be taken to be the image of a torsion unit in
$\ZZ G$ under an ordinary or modular representation of $G$.)
If $F=k$ we assume that $n$ is not divisible by the characteristic $p$ of $k$.
We also assume that $K$ contains a primitive $n$-th root of unity $\zeta$.
The canonical map $R\rightarrow R/\pi R=k$ will be denoted by
$r\mapsto\bar{r}$ ($r\in R$); if $F=k$, it induces an isomorphism
$\lara{\zeta}\cong\lara{\bar{\zeta}}$.
The method to be described rests on the observation that
the eigenvalues of $U$ including their multiplicities can be computed
by Fourier inversion from the traces of the powers of $U$.
Let $\xi$ be an $n$-th root of unity in $K$.
Let $\mu_{\xi}$ be the multiplicity of
$\xi$ as an eigenvalue of $U$ or $\mu_{\bar{\xi}}$ be the multiplicity of
$\bar{\xi}$ as an eigenvalue of $U$, as the case may be.
Having the possibility to diagonalize $U$ shows that
\[ \mu_{\xi} = \frac{1}{n}\sum_{i=0}^{n}\text{\rm trace}(U^{i})\xi^{-i} 
\quad\text{or}\quad
\mu_{\bar{\xi}} = \frac{1}{n}\sum_{i=0}^{n}\lambda(U^{i})\xi^{-i} \]
where in the first case the ordinary trace of the matrices is taken
while in the second $\lambda(U^{i})$ is, in the sense of Brauer, 
the sum of the roots of unity in $\lara{\zeta}$ corresponding 
to all the eigenvalues of $U^{i}$.
Collecting summands for which the matrices $U^{i}$ have the same order,
i.e., taking the Galois action into account yields
\begin{align*}
\mu_{\xi} & = \frac{1}{n}\sum_{d\mid n}
\sum_{\sigma\in\text{\rm Gal}(\sQQ(\zeta^{d})/\sQQ)}
\text{\rm trace}(U^{d})^{\sigma}(\xi^{-d})^{\sigma} \\
& = \frac{1}{n}\sum_{d\mid n}
\text{\rm Tr}_{\sQQ(\zeta^{d})/\sQQ}(\text{\rm trace}(U^{d})\xi^{-d}) \\
\intertext{if $U$ is written over $K$, and in the modular case} 
\mu_{\bar{\xi}} & = \frac{1}{n}\sum_{d\mid n}
\text{\rm Tr}_{\sQQ(\zeta^{d})/\sQQ}(\lambda(U^{d})\xi^{-d}).
\end{align*}

Now let $u$ be a torsion unit of order $n$ in $\NU{\ZZ G}$. Let $\chi$ be an 
ordinary character of $G$ afforded by a representation $D$ (written over
$K$), and let $\varphi$ be a Brauer character of $G$ afforded by a
$p$-modular representation $\Theta$ (written over $k$).
We write $\mu(\xi,u,\chi)$ for the multiplicity of $\xi$ as an eigenvalue
of $D(u)$, and $\mu(\xi,u,\varphi)$ for the multiplicity of $\bar{\xi}$
as an eigenvalue of $\Theta(u)$. Let $\psi$ stand for $\chi$ or $\varphi$.
If $\varphi$ is considered we assume that $p$ does not divide the order
$n$ of $u$, but the order of $G$ (the latter to avoid triviality).
We have seen that
\[ \mu(\xi,u,\psi) = \frac{1}{n}\sum_{d\mid n}
\text{\rm Tr}_{\sQQ(\zeta^{d})/\sQQ}(\psi(u^{d})\xi^{-d}). \]
We split the sum into the two parts
\begin{align*}
a(\xi,u,\psi) & = \frac{1}{n}\sum_{\substack{d\mid n \\ d\neq 1}}
\text{\rm Tr}_{\sQQ(\zeta^{d})/\sQQ}(\psi(u^{d})\xi^{-d}), \\
b(\xi,u,\psi) & = \frac{1}{n}
\text{\rm Tr}_{\sQQ(\zeta)/\sQQ}(\psi(u)\xi^{-1}).
\end{align*}
Note that the multiples $n\cdot a(\xi,u,\psi)$ and $n\cdot b(\xi,u,\psi)$ 
are rational integers.

We immediately note the following useful relation 
for units of prime power order:
\begin{align}\label{Eq4a}
a(\xi,u,\psi)=\frac{1}{r}\mu(\xi^{r},u^{r},\psi) \quad
\text{if $u$ has order a power of a prime $r$.}
\end{align}

If one tries to verify the Zassenhaus conjecture for $G$, i.e., if one 
tries to show that $u$ is rationally conjugate to a group element, then one
can assume---by induction on the order of $u$---that the values
$a(\xi,u,\psi)$ are ``known.'' 
Moreover, since the partial augmentations of $u$ are rational numbers,
we obviously have
\[ b(\xi,u,\chi)=\sum_{x^{G}}\varepsilon_{x}(u)\,b(\xi,x,\chi) \]
where the sum runs over the conjugacy classes of $G$.
Further, it follows from Theorem~\ref{MT} that
\[ b(\xi,u,\varphi)=
\sum_{\substack{x^{G}:\; x\text{ is}\\ \text{$p$-regular}}}
\varepsilon_{x}(u)\,b(\xi,x,\varphi). \]
Altogether, we obtain the equation given by Luthar and Passi:
\begin{align}\label{Eq4}
\mu(\xi,u,\chi) = a(\xi,u,\chi)+\sum_{x^{G}}\varepsilon_{x}(u)\,b(\xi,x,\chi) 
\end{align}
and its modular counterpart
\begin{align}\label{Eq5}
\mu(\xi,u,\varphi) = a(\xi,u,\varphi)+
\sum_{\substack{x^{G}:\; x\text{ is}\\ \text{$p$-regular}}}
\varepsilon_{x}(u)\,b(\xi,x,\varphi). 
\end{align}
The values $b(\xi,x,\chi)$ and $b(\xi,x,\varphi)$ 
can be computed from the (ordinary) character table of $G$ and the 
Brauer character table of $G$ (modulo $p$), respectively.
Both (\ref{Eq4}) and (\ref{Eq5}) should be seen as a 
linear system of equations, indexed by the irreducible characters,
in the unknown partial augmentations $\varepsilon_{x}(u)$. 
Note that the natural integers $\mu(\xi,u,\chi)$ are bounded above by 
the degree of $\chi$.

If $G$ is $p$-solvable, (\ref{Eq5}) provides nothing new since then each 
irreducible Brauer character is the
``reduction modulo $p$'' of an ordinary irreducible character, by the
Fong--Swan--Rukolaine Theorem (see \cite[(22.1)]{CuRe:81}).
For other groups, (\ref{Eq5}) may impose 
stronger conditions on the partial augmentations than (\ref{Eq4}).
This is illustrated by means of a few examples in the next sections.

\section{The symmetric group $S_{5}$}\label{Sec:ex1}

Luthar and Trama \cite{LuTr:91} verified the Zassenhaus conjecture for the 
symmetric group $S_{5}$. Their proof is based on a computer
calculation which yields certain congruences between matrix entries 
with reference to a particular Wedderburn embedding of $\ZZ S_{5}$, 
and uses a few ad hoc arguments.
Here, we show that Theorem~\ref{MT} provides enough information to 
yield another proof.

For torsion units in $\NU{\ZZ S_{5}}$ which are not of order $2$, $4$ or $6$,
the Luthar--Passi method (\ref{Eq4}) yields rational conjugacy to group elements
(as noted in \cite{LuTr:91}), so we only deal with the remaining cases,
for which we will make use only of the Brauer character table for
the prime $5$, shown in Table~\ref{Table2} together with the
decomposition matrix, in the form obtained by requiring 
{\tt CharacterTable("S5") mod 5} in GAP \cite{GAP4} (dots indicate zeros).
\begin{table}[h] 
\[ 
\begin{array}[b]{r|r@{\hspace*{5pt}}r@{\hspace*{5pt}}r@{\hspace*{5pt}}
r@{\hspace*{5pt}}r@{\hspace*{5pt}}r} \hline
& \varphi_{1a} & \varphi_{1b} & \varphi_{3a} & \varphi_{3b}
& \varphi_{5a} & \varphi_{5b} \rule[-7pt]{0pt}{20pt} \\ \hline
\chi_{1a} & 1 & . & . & . & . & . \rule[0pt]{0pt}{13pt} \\
\chi_{1b} & . & 1 & . & . & . & . \\
\chi_{6} & . & . & 1 & 1 & . & . \\
\chi_{4a} & 1 & . & 1 & . & . & . \\
\chi_{4b} & . & 1 & . & 1 & . & . \\
\chi_{5a} & . & . & . & . & 1 & . \\
\chi_{5b} & . & . & . & . & . & 1 \rule[-7pt]{0pt}{5pt} \\
\hline
\end{array}
\qquad
\begin{array}[b]{r|r@{\hspace*{5pt}}r@{\hspace*{5pt}}r@{\hspace*{5pt}}
r@{\hspace*{5pt}}r@{\hspace*{5pt}}r@{\hspace*{5pt}}r} \hline
& {\tt 1a} & {\tt 2a} & {\tt 3a} & {\tt 2b} & {\tt 4a} & {\tt 6a}
\rule[-7pt]{0pt}{20pt} \\ \hline
\varphi_{1a} & 1 & 1 & 1 & 1 & 1 & 1 \rule[0pt]{0pt}{13pt} \\
\varphi_{1b} & 1 & 1 & 1 & -1 & -1 & -1 \\
\varphi_{3a} & 3 & -1 & . & 1 & -1 & -2 \\
\varphi_{3b} & 3 & -1 & . & -1 & 1 & 2 \\
\varphi_{5a} & 5 & 1 & -1 & 1 & -1 & 1 \\
\varphi_{5b} & 5 & 1 & -1 & -1 & 1 & -1 \rule[-7pt]{0pt}{5pt} \\
\hline
\end{array} 
\] 
\rule[0pt]{0pt}{5pt}
\vspace*{-6pt}
\caption{Decomposition matrix and Brauer character table 
of $S_{5}$ ($\text{\rm mod\ }5$)}
\label{Table2}
\vspace*{-10pt}
\end{table}
Let $u$ be a torsion unit in $\NU{\ZZ S_{5}}$,
of order $2$, $4$ or $6$. Let 
$\varepsilon_{{\tt 1a}}, \varepsilon_{{\tt 2a}},
\varepsilon_{{\tt 2b}}, \varepsilon_{{\tt 4a}},
\varepsilon_{{\tt 3a}}, \varepsilon_{{\tt 5a}}$
and $\varepsilon_{{\tt 6a}}$ be its partial augmentations,
so that $\varepsilon_{{\tt 2a}}$, for example, 
denotes the partial augmentation of $u$ with respect to one of the two 
conjugacy classes of elements of order $2$. We have 
$\varepsilon_{{\tt 1a}}=0$ by the Berman--Higman result.
By \cite[Theorem~2.7]{MaRiSeWe:87} (cf.\ Theorem~\ref{T8}), we have 
$\varepsilon_{{\tt 3a}}=\varepsilon_{{\tt 5a}}=\varepsilon_{{\tt 6a}}=0$
if $u$ is of order $2$ or $4$, and $\varepsilon_{{\tt 5a}}=0$ if $u$ is of
order $6$.

We denote by $\Theta_{\ast}$ a representation of $G$, over a sufficiently
large field, which affords the Brauer character $\varphi_{\ast}$, and
write $\Theta_{\ast}(u)\sim\text{\rm diag}(\lambda_{1},\dotsc,
\lambda_{\varphi_{\ast}(1)})$ to indicate that $\Theta_{\ast}(u)$ 
has eigenvalues $\lambda_{1},\dotsc,\lambda_{\varphi_{\ast}(1)}$ 
(multiplicities taken into account), which will, as happened before, 
be identified with complex roots of unity.

\subsection*{When $\boldsymbol{u}$ has order $\boldsymbol{2}$.}
To simplify things, we shall use that $\varepsilon_{{\tt 4a}}=0$ for 
theoretical reason, viz.\ Theorem~\ref{T8}. 
Then taking augmentation of $u$ gives 
$\varepsilon_{{\tt 2a}}+\varepsilon_{{\tt 2b}}=1$, and
$\varphi_{1b}(u)=\varepsilon_{{\tt 2a}}-\varepsilon_{{\tt 2b}}=
2\varepsilon_{{\tt 2a}}-1$. Note that $\varphi_{1b}(u)=\pm 1$ since 
$\varphi_{1b}$ is of degree one. It follows that 
$\varepsilon_{{\tt 2a}}=1$ or $\varepsilon_{{\tt 2a}}=0$, and we are done.

\subsection*{When $\boldsymbol{u}$ has order $\boldsymbol{4}$.}
Taking augmentation of $u$ gives 
\begin{align}\label{E1}
\varepsilon_{{\tt 2a}}+\varepsilon_{{\tt 2b}}+\varepsilon_{{\tt 4a}}=1.
\end{align}
Since $u^{2}$ is rationally conjugate to ${\tt 2a}$ or ${\tt 2b}$,
the matrix $\Theta_{{\tt 3a}}(u^{2})$ is conjugate to 
$\Theta_{{\tt 3a}}({\tt 2a})$ or $\Theta_{{\tt 3a}}({\tt 2b})$
by Remark~\ref{remi2}, and as
both $\varphi_{{\tt 3a}}({\tt 2a})$ and $\varphi_{{\tt 3a}}({\tt 2b})$
are not equal to the degree of $\varphi_{{\tt 3a}}$, it follows that
$\Theta_{{\tt 3a}}(u)$ is of order $4$. Together with 
$\Theta_{{\tt 3a}}(u)\in\ZZ$, this implies  
$\Theta_{{\tt 3a}}(u)\sim\text{\rm diag}(i,-i,\pm 1)$, so
$\varphi_{{\tt 3a}}(u)=\pm 1$. Similarly, 
$\varphi_{{\tt 3b}}(u)=\pm 1$. Thus
\begin{align}
\varphi_{{\tt 3a}}(u) & = 
-\varepsilon_{{\tt 2a}}+\varepsilon_{{\tt 2b}}-\varepsilon_{{\tt 4a}}=\pm 1, 
\label{E3}
\\
\label{E4}
\varphi_{{\tt 3b}}(u) & = 
-\varepsilon_{{\tt 2a}}-\varepsilon_{{\tt 2b}}+\varepsilon_{{\tt 4a}}=\pm 1. 
\end{align}
Adding (\ref{E1}) and (\ref{E3}) shows that
\begin{align}\label{E2}
\varepsilon_{{\tt 2b}}=1 \quad \text{or \quad}
\varepsilon_{{\tt 2b}}=0. 
\end{align}
If $\varphi_{{\tt 3b}}(u)=-1$, then adding (\ref{E1}) and (\ref{E4}) 
shows that $\varepsilon_{{\tt 4a}}=0$, which
contradicts \cite[(7.3)]{Seh:93}. Thus $\varphi_{{\tt 3b}}(u)=1$, 
and adding (\ref{E1}) and (\ref{E4}) shows that $\varepsilon_{{\tt 4a}}=1$.
Consequently $\varepsilon_{{\tt 2a}}+\varepsilon_{{\tt 2b}}=0$ by 
(\ref{E1}). So $\varepsilon_{{\tt 2a}}=-1$ if $\varepsilon_{{\tt 2b}}=1$. But
then $\varphi_{{\tt 1b}}(u)=
\varepsilon_{{\tt 2a}}-\varepsilon_{{\tt 2b}}-\varepsilon_{{\tt 4a}}=
-1-1-1=-3$, contradicting  $\varphi_{{\tt 1b}}(u)=\pm 1$.
Hence $\varepsilon_{{\tt 2b}}=0$ by (\ref{E2}), and 
$\varepsilon_{{\tt 2a}}=0$, so $\varepsilon_{{\tt 4a}}$ is the only
partial augmentation which does not vanish.

\subsection*{When $\boldsymbol{u}$ has order $\boldsymbol{6}$.}
Taking augmentation of $u$ gives 
\begin{align}\label{E5}
\varepsilon_{{\tt 2a}}+\varepsilon_{{\tt 2b}}+\varepsilon_{{\tt 4a}}
+\varepsilon_{{\tt 3a}}+\varepsilon_{{\tt 6a}}=1.
\end{align}
Also,
\begin{align}\label{E6}
\varphi_{{\tt 1b}}(u)=
\varepsilon_{{\tt 2a}}-\varepsilon_{{\tt 2b}}-\varepsilon_{{\tt 4a}}
+\varepsilon_{{\tt 3a}}-\varepsilon_{{\tt 6a}}=\pm 1.
\end{align}
The unit $u^{3}$ of order $2$ is rationally conjugate to either ${\tt 2a}$
or ${\tt 2b}$. If $u^{3}$ is conjugate to ${\tt 2a}$, then $u^{3}$, and
hence also $u$, maps to $1$ under the natural map 
$\ZZ S_{5}\rightarrow \ZZ S_{5}/A_{5}=\ZZ C_{2}$, so
$\varphi_{{\tt 1b}}(u)=1$. Adding (\ref{E5}) and (\ref{E6}) then shows:
\begin{align}\label{E7}
\text{$u^{3}$ is conjugate to ${\tt 2a}$}\;\,\Longrightarrow \;\,
\varepsilon_{{\tt 2a}}+\varepsilon_{{\tt 3a}}=1.
\end{align}
Also, if $u^{3}$ is conjugate to ${\tt 2b}$, then $u^{3}$, and
hence also $u$, maps to $-1$ under the natural map 
$\ZZ S_{5}\rightarrow \ZZ C_{2}$, so $\varphi_{{\tt 1b}}(u)=-1$, and 
adding (\ref{E5}) and (\ref{E6}) shows:
\begin{align}\label{E8}
\text{$u^{3}$ is conjugate to ${\tt 2b}$}\;\,\Longrightarrow \;\,
\varepsilon_{{\tt 2a}}+\varepsilon_{{\tt 3a}}=0.
\end{align}

Let $\Theta$ be one of $\Theta_{{\tt 3a}}$ and $\Theta_{{\tt 3b}}$, with
Brauer character $\varphi$. Then we can write
\[ \Theta(u)\sim\text{\rm diag}(\nu_{1}\zeta^{\alpha_{1}},
\nu_{2}\zeta^{\alpha_{2}},\nu_{3}\zeta^{\alpha_{3}}) \] for some 
$\nu_{i}\in\{\pm 1\}$, a primitive third root of unity $\zeta$, and
$\alpha_{i}\in\ZZ$. Since $\varphi(u^{2})=\varphi({\tt 3a})=0$
we have $\zeta^{\alpha_{1}}+\zeta^{\alpha_{2}}+\zeta^{\alpha_{3}}=0$.
Together with $\Theta(u)\in\ZZ$, it follows that
\begin{align}\label{E9}
\Theta(u)\sim\pm\text{\rm diag}(1,-\zeta,-\zeta^{2}). 
\end{align}
If $\Theta(u)\sim\text{\rm diag}(-1,\zeta,\zeta^{2})$, then 
$\varphi(u^{3})=1$, so $\varphi=\varphi_{{\tt 3a}}$ and $u^{3}$ 
is conjugate to ${\tt 2b}$.
Suppose that $u^{3}$ is conjugate to ${\tt 2a}$.
Then the preceding observations show that 
$\varphi_{{\tt 3a}}(u)=\varphi_{{\tt 3b}}(u)=2$.
It follows that $-2\varepsilon_{{\tt 2a}}=\varphi_{{\tt 3a}}(u)+
\varphi_{{\tt 3b}}(u)=4$, and $\varepsilon_{{\tt 2a}}=-2$.
Now $\varepsilon_{{\tt 3a}}=3$ from (\ref{E7}). But then 
$\varphi_{{\tt 5a}}(u)+\varphi_{{\tt 5b}}(u)=
2\varepsilon_{{\tt 2a}}-2\varepsilon_{{\tt 3a}}=-4-6=-10$, 
implying $\varphi_{{\tt 5a}}(u^{2})=\varphi_{{\tt 5a}}(1)$
which is impossible.

Hence $u^{3}$ is conjugate to ${\tt 2b}$. So
$\varphi_{{\tt 3a}}(u^{3})=1$ and $\varphi_{{\tt 3b}}(u^{3})=-1$,
and (\ref{E9}) gives $\varphi_{{\tt 3a}}(u)=-2$ and
$\varphi_{{\tt 3b}}(u)=2$.
It follows that $-2\varepsilon_{{\tt 2a}}=\varphi_{{\tt 3a}}(u)+
\varphi_{{\tt 3b}}(u)=0$, i.e., $\varepsilon_{{\tt 2a}}=0$.
Now $\varepsilon_{{\tt 3a}}=0$ from (\ref{E8}).
Using (\ref{E5}), which now reads
$\varepsilon_{{\tt 2b}}+\varepsilon_{{\tt 4a}}
+\varepsilon_{{\tt 6a}}=1$, we get
$\varphi_{{\tt 5a}}(u)=\varepsilon_{{\tt 2b}}-\varepsilon_{{\tt 4a}}
+\varepsilon_{{\tt 6a}}=1-2\varepsilon_{{\tt 4a}}$.
Subtracting (\ref{E5}) from 
$\varepsilon_{{\tt 2b}}-\varepsilon_{{\tt 4a}}
-2\varepsilon_{{\tt 6a}}=\varphi_{{\tt 3a}}(u)=-2$ gives
$2\varepsilon_{{\tt 4a}}+3\varepsilon_{{\tt 6a}}=3$, and we further obtain
$\varphi_{{\tt 5a}}(u)=-2+3\varepsilon_{{\tt 6a}}$.
Since $\abs{-2+3\varepsilon_{{\tt 6a}}}<\varphi_{{\tt 5a}}(1)=5$, it follows
that $\varepsilon_{{\tt 6a}}\in\{0,1,2\}$. If 
$\varepsilon_{{\tt 6a}}=0$ or $\varepsilon_{{\tt 6a}}=2$, then
$\varepsilon_{{\tt 4a}}=3/2$ or $\varepsilon_{{\tt 4a}}=-3/2$, which is
impossible. Hence $\varepsilon_{{\tt 6a}}=1$, and  
$\varepsilon_{{\tt 4a}}=0$, $\varepsilon_{{\tt 2b}}=0$ easily follows.
Thus $\varepsilon_{{\tt 6a}}$ is the only
partial augmentation which does not vanish.

\begin{remark}
Suppose that the Zassenhaus conjecture is verified for a symmetric group
$S_{n}$. Then the Zassenhaus conjecture also holds for the 
symmetric groups $S_{m}$, $m<n$, since there is no fusion of conjugacy
classes of $S_{m}$ in $S_{n}$. 
\end{remark}

\section{The groups $\text{\rm PSL}(2,p^{f})$}\label{Sec:ex2}

The projective special linear group $\text{\rm PSL}(2,p^{f})$ 
may be seen as a prototype for intended applications. 
Most of its conjugacy classes are $p$-regular, and a $p$-singular group
element is of order $p$, so it seems particular promising
to work with Brauer characters in the defining characteristic $p$.
These characters are known and it is also advantageous that
there exist some of small degree.

We remark that $\text{\rm PSL}(2,3)$ and $\text{\rm PSL}(2,5)$ are the 
alternating groups of order $12$ and $60$, respectively, for which
the Zassenhaus conjecture has been verified \cite{LuPa:89}. 

Group-theoretical properties of $\text{\rm PSL}(2,p^{f})$ are described in
\cite[Kapitel~II, \S8]{Hup:67}.
Its order is $\frac{1}{k}(p^{f}-1)p^{f}(p^{f}+1)$ where
$k=\text{\rm gcd}(p^{f}-1,2)$. Of importance to us will be that
$\text{\rm PSL}(2,p^{f})$ has cyclic subgroups of order
$\frac{1}{k}(p^{f}-1)$ and $\frac{1}{k}(p^{f}+1)$, and that each $p$-regular 
element of $\text{\rm PSL}(2,p^{f})$ lies in a conjugate of one of 
these subgroups. Moreover, if $g$ is a $p$-regular element of 
$\text{\rm PSL}(2,p^{f})$ of order greater than $2$, then
$g^{-1}$ is its only distinct conjugate in $\lara{g}$.

The representation theory of $\text{\rm PSL}(2,p^{f})$ is well
understood, and some references are given below (we will only use characters).
The $p$-modular representation theory of 
$\text{\rm SL}(2,p)$ is described in Alperin's beautiful book \cite{Alp:86}, 
and there is also a nice introductory article \cite{Hum:75} by Humphreys 
on the representations of $\text{\rm SL}(2,p)$.
The character table of $\text{\rm SL}(2,p)$ was first obtained by Frobenius;
shortly afterwards, Schur and (independently) H.~Jordan found the
characters of $\text{\rm SL}(2,p^{f})$.
Dornhoff's version \cite[\S38]{Dor:71} is a readable account of Schur's work.
The modular irreducible representations of $\text{\rm SL}(2,p^{f})$
in describing characteristic were given by Brauer and Nesbitt
\cite[\S30]{BrNe:41}; they also gave the decomposition matrix for 
$\text{\rm PSL}(2,p)$, $p$ odd. 
The decomposition matrices for all groups
$\text{\rm PSL}(2,p^{f})$ were given by Burkhardt \cite{Bur:76}
(see also \cite{HSW:82}).
In \cite{Sri:64}, the Brauer character table (relative to $p$) of 
$\text{\rm SL}(2,p^{f})$ (of which the Brauer character table of
$\text{\rm PSL}(2,p^{f})$ is a part of) is given explicitly.

First, we shall consider $p$-singular torsion units when $p$ is odd. 
For $\text{\rm PSL}(2,p)$ one can give definite results.
The ordinary character table of $\text{\rm PSL}(2,p^{f})$ for odd $p$
is shown in Table~\ref{Table3}.
\newlength{\uni}
\settowidth{\uni}{$- $}
\newlength{\unii}
\settowidth{\unii}{$\varepsilon=-1:\; {}$}
\begin{table}[h] 
\[ \begin{array}{c}
\begin{array}{c|cccccc} \hline
\text{class of} & 1 & c & d & a^{l} & b^{m} \rule[-7pt]{0pt}{20pt} \\ \hline
\text{order} 
& 1 & p & p & 
\text{\rm o}(a)=
\frac{q-1}{2} & \text{\rm o}(b)=\frac{q+1}{2}
\rule[-7pt]{0pt}{20pt} \\ \hline
1 & 1 & 1 & 1 & 1 & 1 \rule[0pt]{0pt}{13pt} \\
\psi & q & 0 & 0 & 1 & -1\hspace*{\uni} \\
\chi_{i} & q+1 & 1 & 1 & \rho^{il}+ \rho^{-il} & 0 \\
\theta_{j} & q-1 & -1\hspace*{\uni} & -1\hspace*{\uni} 
& 0 & -(\sigma^{jm}+ \sigma^{-jm}) \\
\eta_{1} & \frac{1}{2}(q+\varepsilon) & 
\frac{1}{2}(\varepsilon+\sqrt{\varepsilon q})
 & \frac{1}{2}(\varepsilon-\sqrt{\varepsilon q})
 & (-1)^{l}\delta_{\varepsilon,1} & (-1)^{m}\delta_{\varepsilon,-1} 
\rule[-7pt]{0pt}{20pt} \\
\eta_{2} & \frac{1}{2}(q+\varepsilon) & 
\frac{1}{2}(\varepsilon-\sqrt{\varepsilon q})
 & \frac{1}{2}(\varepsilon+\sqrt{\varepsilon q})
& (-1)^{l}\delta_{\varepsilon,1} & (-1)^{m}\delta_{\varepsilon,-1} 
\rule[-7pt]{0pt}{5pt} \\ \hline
\end{array} \\
\begin{array}{rl}
\text{Entries:}\ &
\text{sign } \varepsilon \text{ such that } q\equiv\varepsilon\mod{4}, \\
& \delta_{\varepsilon,\pm 1} \text{ Kronecker symbol,} \\ 
& \rho=e^{4\pi i/(q-1)},\;\sigma=e^{4\pi i/(q+1)}, \\
& \parbox[t]{\unii}{$\varepsilon=1:\; $ \hfill} 
1\leq i\leq \frac{1}{4}(q-5),\; 1\leq j,l,m\leq \frac{1}{4}(q-1),
\rule[6pt]{0pt}{5pt} \\
& \varepsilon=-1:\;
1\leq i,j,l\leq \frac{1}{4}(q-3),\; 1\leq m\leq \frac{1}{4}(q+1). 
\rule[6pt]{0pt}{5pt} \\
\end{array}
\rule[0pt]{0pt}{50pt}
\end{array} \] 
\rule[0pt]{0pt}{5pt}
\caption{Character table of $\text{\rm PSL}(2,q)$, 
$q=p^{f}\geq 5$, odd prime $p$}
\label{Table3}
\vspace*{-10pt}
\end{table}

The following two results were proved by Wagner \cite{Wag:95} 
who did the formal calculations needed for application of 
the Luthar--Passi method (\ref{Eq4}). The calculation 
of (\ref{Eq4}) for the following result in the case $f=1$ is reported 
in \cite{BHK:04}. We shall give more direct arguments.
\begin{proposition}\label{PSLa}
Let $G=\text{\rm PSL}(2,p^{f})$ for an odd prime $p$, and $f\leq 2$.
Then units of order $p$ in $\NU{\ZZ G}$ are rationally conjugate to
elements of $G$.
\end{proposition}
\begin{proof}
Let $u$ be a torsion unit in $\NU{\ZZ G}$ of order $p$, and let $\alpha$
be its partial augmentation with respect to the class of $c$, an element 
of order $p$. Then $1-\alpha$ is its partial augmentation with respect 
to the other class of elements of order $p$.

First, suppose that $G=\text{\rm PSL}(2,p)$, and set $\zeta=e^{2\pi i/p}$.
The character table of $G$ involves Gaussian sums:
\[ \tfrac{1}{2}(\varepsilon+\sqrt{\varepsilon p})=\delta_{\varepsilon,1}+\!
\sum_{\substack{\text{$i$ square}\\ \text{in $(\ZZ/p\ZZ)^{\times}$}}}
\!\zeta^{i}, \qquad
\tfrac{1}{2}(\varepsilon-\sqrt{\varepsilon p})=\delta_{\varepsilon,1}+\!
\sum_{\substack{\text{$j$ non-square}\\ \text{in $(\ZZ/p\ZZ)^{\times}$}}}
\!\zeta^{j} \]
(see, for example, \cite[Chapter~VI, Theorem~3.3]{Lan:02}). Thus 
\begin{align*}
\eta_{1}(u) & =\alpha\tfrac{1}{2}(\varepsilon+\sqrt{\varepsilon p})+
(1-\alpha)\tfrac{1}{2}(\varepsilon-\sqrt{\varepsilon p}) \\
& = \delta_{\varepsilon,1}+(\alpha-1)+(2\alpha-1)
\sum_{\substack{\text{$i$ square}\\ \text{in $(\ZZ/p\ZZ)^{\times}$}}}
\!\zeta^{i}.
\end{align*}
If $\alpha>0$, then this character value is a sum of 
$\tfrac{1}{2}(p+\varepsilon)$ roots of unity if and only if 
$\alpha=1$.
By symmetry (we may also use $\eta_{2}$), it follows
that one of the partial augmentations on elements of order $p$ vanish.

Now suppose that $G=\text{\rm PSL}(2,p^{2})$. Then $\varepsilon=1$ as
$p^{2}\equiv 1\mod{4}$, and we have
\begin{align*}
\eta_{1}(u) & =\alpha\tfrac{1}{2}(1+p)+
(1-\alpha)\tfrac{1}{2}(1-p) \\
& = \tfrac{1}{2}((2\alpha-1)p+1).
\end{align*}
Note that if $n$ is a natural number, and $-n$ is a sum of $m$ 
$p$-th roots of unity, then $m\geq n(p-1)$. Hence if $\alpha<0$, then
\begin{align*}
\tfrac{1}{2}(p^{2}+1) & =\eta_{1}(1)\geq
\abs{\eta_{1}(u)}(p-1) = \tfrac{1}{2}(-(2\alpha-1)p-1)(p-1) \\
& \geq \tfrac{1}{2}(3p-1)(p-1) 
> \tfrac{1}{2}(p^{2}+1),
\end{align*}
which is impossible. Thus $\alpha\geq 0$. By symmetry (use $\eta_{2}$),
also $1-\alpha\geq 0$. So $\alpha=0$ or $\alpha=1$, and we are done.
\end{proof}

\begin{remark}
Let $G=\text{\rm PSL}(2,p^{f})$ for an odd prime $p$, and set
$H=\text{\rm PSL}(2,p^{2f})$. Considering $G$ as a subgroup of $H$,
the two conjugacy classes of elements of order $p$ in $G$ fuse
in $H$. Thus a unit of order $p$ in $\NU{\ZZ G}$ is conjugate to
an element of $G$ by a unit of $\QQ H$.
\end{remark}

Apparently, we exploited the fact that some character degrees are
small when compared with the orders of the involved elements. 
This is also what underlies the proof of the next proposition.
\begin{proposition}\label{PSLb}
Let $G=\text{\rm PSL}(2,p)$ for an odd prime $p$. If the order of a torsion
unit $u$ in $\NU{\ZZ G}$ is divisible by $p$, then $u$ is of order $p$.
\end{proposition}
\begin{proof}
We can assume that $p>5$ (see \cite{LuPa:89}). 

Since $G$ has no elements of order $p^{2}$, there are no 
torsion units of order $p^{2}$ in $\NU{\ZZ G}$ (see \cite[(7.3)]{Seh:93}).
Suppose that a torsion unit $u$ in $\NU{\ZZ G}$ is of order $pr$, for
some prime $r$ different from $p$. Let $\varepsilon_{c}$ and 
$\varepsilon_{d}$ be the partial augmentations of $u$ with respect to
the two conjugacy classes of $G$ of elements of order $p$. Write
$u=\sum_{g\in G}a_{g}g$ (all $a_{g}$ in $\ZZ$).
Elementary calculus shows (see \cite[\S7]{Seh:93}):
\[ u^{p}\in \sum_{g\in G}a_{g}g^{p} + [\ZZ G,\ZZ G] + p\ZZ G. \]
Thus the $1$-coefficient of $u^{p}$ lies in 
$\varepsilon_{c}+\varepsilon_{d}+p\ZZ$. By the Berman--Higman result, 
it follows that $\varepsilon_{c}+\varepsilon_{d}=mp$ for some integer $m$. 
Similarly, we obtain that the sum of all partial augmentations of $u$ 
with respect to conjugacy classes of elements of order $r$ is divisible 
by $r$. By Theorem~\ref{T8}, all other partial augmentations of $u$ vanish.
Thus $m\neq 0$ as $u$ has augmentation one.
Suppose that $r\mid \frac{1}{2}(p-1)$. Then $\theta_{1}(x)=0$ for all 
$r$-elements $x$ of $G\setminus\{1\}$, so $\theta_{1}(u)=-mp$, which is 
impossible since $\theta_{1}$ has degree $p-1$. Hence 
$r\mid \frac{1}{2}(p+1)$. So $\chi_{1}(x)=0$ for all $r$-elements $x$ of 
$G\setminus\{1\}$, and $\chi_{1}(u)=\varepsilon_{c}+\varepsilon_{d}=mp$. 
As $\chi_{1}$ has degree $p+1$, it follows that
$m=\pm 1$, and $\pm p$ is the sum of $p+1$ ($pr$)-th roots of unity.
This can only happen if $p=3$ (and $r=2$), a case which we excluded
at the beginning. The proposition is proved.
\end{proof}

We now consider $p$-regular torsion units. 
We only prove partial results, but perhaps more can be expected when 
the $p$-modular version of the Luthar--Passi method is used more rigorously.

Let $a$ and $b$ be elements of $\text{\rm PSL}(2,p^{f})$
of order $\frac{1}{k}(p^{f}-1)$ and 
$\frac{1}{k}(p^{f}+1)$, respectively, where $k=\text{\rm gcd}(p^{f}-1,2)$.
Let $\rho$ and $\sigma$ be roots of unity of the same order as 
$a$ and $b$, respectively. Relative to the prime $p$, the group 
$\text{\rm PSL}(2,p^{f})$ has 
Brauer characters $\varphi_{1},\varphi_{3},\varphi_{5},\ldots$ of degree 
$1,3,5,\ldots$ arising from the action of $\text{\rm SL}(2,p^{f})$ on 
homogeneous polynomials in two variables of degree $0,2,4,\ldots$, such 
that if $\Theta_{2l+1}$ is a representation affording $\varphi_{2l+1}$, 
then (in the notation from Section~\ref{Sec:ex1}):
\begin{align*}
\Theta_{2l+1}(a) & \sim\text{\rm diag}(\rho^{l},\rho^{l-1},\dotsc,
\rho,1,\rho^{-1},\dotsc,\rho^{-(l-1)},\rho^{-l}), \\
\Theta_{2l+1}(b) & \sim\text{\rm diag}(
\sigma^{l},\sigma^{l-1},\dotsc,\sigma,1,\sigma^{-1},\dotsc,
\sigma^{-(l-1)},\sigma^{-l}).
\end{align*}

The proof of our first result is straightforward.
\begin{proposition}\label{primeord}
Let $G=\text{\rm PSL}(2,p^{f})$, and let $r$ be a prime different from $p$.
Then any torsion unit $u$ in $\NU{\ZZ G}$ of order $r$ is
rationally conjugate to an element of $G$. 
\end{proposition}
\begin{proof}
Only for elements $x$ of $G$ which are of order $r$ we possibly have 
$\varepsilon_{x}(u)\neq 0$, by Theorem~\ref{T8}.
So we can assume that $r\neq 2$ as there is only one conjugacy class
of elements of order $2$ in $G$.
Let $x$ be an element of order $r$ in $G$. Then 
$\varphi_{3}(x)=1+\zeta+\zeta^{-1}$ for a primitive $r$-th root of unity
$\zeta$. By Theorem~\ref{MT},
\begin{align*}
\varphi_{3}(u) & =
\sum_{i=1}^{(r-1)/2}\varepsilon_{x^{i}}(u)(1+\zeta^{i}+\zeta^{-i}) \\
& = 1+\sum_{i=1}^{(r-1)/2}\varepsilon_{x^{i}}(u)(\zeta^{i}+\zeta^{-i}). 
\end{align*}
In particular, $\varphi_{3}(u)\in\ZZ[\zeta+\zeta^{-1}]$, and since
$\varphi_{3}(u)$ is a sum of three $r$-th roots of unity, 
$\varphi_{3}(u)=1+\zeta^{j}+\zeta^{-j}$ with $0\leq j\leq (r-1)/2$.
Comparing both expressions for $\varphi_{3}(u)$, it follows that 
$j\neq 0$ (otherwise all $\varepsilon_{x^{i}}(u)$ would be equal to $-2$,
which is impossible), and further $\varepsilon_{x^{j}}(u)=1$, 
with all other partial augmentations of $u$ vanishing.
\end{proof}

The next proposition confirms a conjecture of Bovdi 
(see Section~\ref{Sec:PA}) for the groups $\text{\rm PSL}(2,p^{f})$.
Its proof is already a little bit more sophisticated.
\begin{proposition}\label{Bov}
Let $G=\text{\rm PSL}(2,p^{f})$, and let $r$ be a prime different from $p$.
Suppose that $u$ is an $r$-torsion unit in $\NU{\ZZ G}$, of order $r^{n}$
(say). For any integer $m$, let $T_{r^{m}}\subset G$ be a set of 
representatives of the conjugacy classes
of elements of order $r^{m}$ in $G$. Then
\begin{align*}
\sum_{x\in T_{r^{m}}}\varepsilon_{x}(u)=0\quad\text{for all $0\leq m<n$.}
\end{align*}
Moreover, there exist elements of $G$ of the same order as $u$, and for 
such an element $g$, we have $\mu(1,u,\varphi)=\mu(1,g,\varphi)$ for each 
$p$-modular character $\varphi$ of $G$ 
(notation as in Section~\ref{Sec:LPmethod}).
\end{proposition}
\begin{proof}
There exists an element $g$ in $G$ of the same order as $u$
(see \cite[(7.3)]{Seh:93}). We shall use induction, firstly on $n$ and 
secondly on $m$, to prove $\mu(1,u,\varphi)=\mu(1,g,\varphi)$ for each 
$p$-modular character $\varphi$ of $G$, and 
$\sum_{x\in T_{r^{m}}}\varepsilon_{x}(u)=0$ for all $0\leq m<n$.
By Proposition~\ref{primeord}, and the Berman--Higman result, this
holds for $n=1$, and we can assume that $n>1$, $m\geq 1$.
Suppose that the assertion holds for smaller values of $n$ or $m$.

We also write $T_{\leq r^{m}}$ for a set of representatives of the conjugacy
classes of $r$-elements of $G$ of order equal or less than $r^{m}$ in $G$, 
and $T_{>r^{m}}$ for such a set of $r$-elements whose order is 
greater than $r^{m}$, but equal or less than $r^{n}$. (The latter condition 
is just a technical one, which makes sense in view of Theorem~\ref{T8}).
Let $\zeta$ be a primitive $r^{n}$-th root of unity,
and set $\text{\rm Tr}=\text{\rm Tr}_{\sQQ(\zeta)/\sQQ}$.
Note that for any integer $l$,
\[ \text{\rm Tr}(\zeta^{l})=\begin{cases}
r^{n-1}(r-1) & \text{if $\zeta^{l}=1$}, \\
-r^{n-1} & \text{if $\zeta^{l}\neq 1$ and $\zeta^{lr}=1$}, \\
0 & \text{otherwise.} \end{cases} \]
Set $k=2r^{m-1}+1$. Then
\begin{align}\label{eqq1} 
\frac{1}{r^{n}}\text{\rm Tr}(\varphi_{k}(x))=\begin{cases}
(r-1)/r & \text{if $x\in T_{>r^{m}}$}, \\
(r-3)/r & \text{if $x\in T_{r^{m}}$}. \end{cases} 
\end{align}
By (\ref{Eq4a}) and (\ref{Eq5}),
\[ \mu(1,u,\varphi_{k})=\frac{1}{r}\mu(1,u^{r},\varphi_{k})+
\frac{1}{r^{n}}\sum_{x\in T_{\leq r^{n}}}\varepsilon_{x}(u)\,
\text{\rm Tr}(\varphi_{k}(x)). \]
Inductively, we may assume that 
$\mu(1,u^{r},\varphi_{k})=\mu(1,g^{r},\varphi_{k})=1$.
Note that for elements $x$ and $y$ of $\lara{g}$ of the same order, 
$\varphi_{k}(x)$ and $\varphi_{k}(y)$ are algebraically conjugate, so
$\text{\rm Tr}(\varphi_{k}(x))=\text{\rm Tr}(\varphi_{k}(y))$.
Hence for $0\leq l<m$, we obtain inductively
\[ \sum_{x\in T_{r^{l}}}\varepsilon_{x}(u)\,\text{\rm Tr}(\varphi_{k}(x))=
\text{\rm Tr}\big(\varphi_{k}\big(g^{r^{n-l}}\big)\big)
\sum_{x\in T_{r^{l}}}\varepsilon_{x}(u)=0. \]
Together with (\ref{eqq1}), it follows that
\begin{align*}
\mu(1,u,\varphi_{k}) & = \frac{1}{r}+\frac{r-3}{r}
\sum_{x\in T_{r^{m}}}\varepsilon_{x}(u)+\frac{r-1}{r}
\sum_{x\in T_{>r^{m}}}\varepsilon_{x}(u) \\
& = \frac{1}{r}+\frac{r-3}{r}
\sum_{x\in T_{\leq r^{m}}}\varepsilon_{x}(u)+\frac{r-1}{r}
\bigg(1-\sum_{x\in T_{\leq r^{m}}}\varepsilon_{x}(u)\bigg) \\ 
& = 1-\frac{2}{r}\sum_{x\in T_{\leq r^{m}}}\varepsilon_{x}(u).
\end{align*}
Hence $\mu(1,u,\varphi_{k})\equiv 1\mod{2}$ since 
$\sum_{x\in T_{\leq r^{m}}}\varepsilon_{x}(u)$ is divisible by $r$,
see (\ref{CLequi}). On the other hand, 
$\mu(1,u,\varphi_{k})\leq 1$ as $\mu(1,u^{r},\varphi_{k})=1$.
So $\mu(1,u,\varphi_{k})=1$, and
$\sum_{x\in T_{\leq r^{m}}}\varepsilon_{x}(u)=0$ as desired.

It remains to show that $\mu(1,u,\varphi)=\mu(1,g,\varphi)$
for a $p$-modular character $\varphi$ of $G$. As above, we obtain inductively
\begin{align*}
\mu(1,u,\varphi) & = \frac{1}{r}\mu(1,u^{r},\varphi)+
\frac{1}{r^{n}}\sum_{x\in T_{\leq r^{n}}}\varepsilon_{x}(u)\,
\text{\rm Tr}(\varphi(x)) \\
 & = \frac{1}{r}\mu(1,g^{r},\varphi)+
\frac{1}{r^{n}}\sum_{x\in T_{r^{n}}}\varepsilon_{x}(u)\,
\text{\rm Tr}(\varphi(x)) \\
& = \frac{1}{r}\mu(1,g^{r},\varphi)+\frac{1}{r^{n}}
\text{\rm Tr}(\varphi(g))\sum_{x\in T_{r^{n}}}\varepsilon_{x}(u) \\
& = \mu(1,g,\varphi).
\end{align*}
The proposition is proved.
\end{proof}

We may expect that results about $p$-regular torsion units of mixed order 
can be achieved as well.
The following proposition was proved in \cite[Example~3.6]{He:04a} for
$\text{\rm PSL}(2,7)$ by an ad hoc argument. 
The proof turns out to be pretty simple.
\begin{proposition}\label{ord6}
Let $G=\text{\rm PSL}(2,p^{f})$ with $p\neq 2,3$. 
Then any torsion unit in $\NU{\ZZ G}$ 
of order $6$ is rationally conjugate to an element of $G$. 
\end{proposition}
\begin{proof}
Suppose that there is a torsion unit $u$ of order $6$ in $\NU{\ZZ G}$.
Then $G$ has elements of order $3$ which form a single conjugacy class.
Of course, the elements of order $2$ are all conjugate in $G$.
Further, if $G$ has elements of order $6$, they also form
a single conjugacy class.
We denote the corresponding partial augmentations of $u$ by
$\varepsilon_{2}$, $\varepsilon_{3}$ and $\varepsilon_{6}$, agreeing
that $\varepsilon_{6}=0$ if $G$ contains no element of order $6$.
By Theorem~\ref{T8}, these are the only partial augmentations of $u$
which are possibly nonzero. 
The Brauer characters $\varphi_{3}$ and $\varphi_{5}$ have the following
values on the relevant classes (the last column may be left out):
\[ \begin{array}{c|rrr} \hline
 & {\tt 2a} & {\tt 3a} & {\tt 6a} \rule[-7pt]{0pt}{20pt} \\ \hline
\varphi_{3} & -1 & 0 & 2 \rule[0pt]{0pt}{13pt} \\
\varphi_{5} & 1 & -1 &  1
\rule[-7pt]{0pt}{5pt} \\ \hline
\end{array}  \]
By Proposition~\ref{primeord} and Remark~\ref{remi2}, we have
\[ \Theta_{5}(u^{3})\sim\text{\rm diag}(1,-1,1,-1,1), \quad
\Theta_{5}(u^{2})\sim
\text{\rm diag}(\zeta^{2},\zeta,1,\zeta^{-1},\zeta^{-2}) \]
for a primitive third root of unity $\zeta$. By Theorem~\ref{MT}, 
$\varphi_{5}(u)=\varepsilon_{2}-\varepsilon_{3}+\varepsilon_{6}$.
In particular $\varphi_{5}(u)\in\ZZ$, forcing 
$\Theta_{5}(u)\sim\text{\rm diag}(\zeta,\zeta^{2},1,-\zeta,-\zeta^{2})$.
Thus $\varphi_{5}(u)=1$. Taking augmentation of $u$ gives
$\varepsilon_{2}+\varepsilon_{3}+\varepsilon_{6}=1$, so that
$\varepsilon_{3}=0$ follows. 
For the same reasons as before, we further have 
\[ \Theta_{3}(u^{3})\sim\text{\rm diag}(-1,1,-1), \quad
\Theta_{3}(u^{2})\sim\text{\rm diag}(\zeta,1,\zeta^{-1}) \]
and $\varphi_{3}(u)=-\varepsilon_{2}+2\varepsilon_{6}$.
This forces $\Theta_{3}(u)\sim\text{\rm diag}(-\zeta^{2},1,-\zeta)$ 
as $\varphi_{3}(u)\in\ZZ$, and $\varphi_{3}(u)=2$. 
Using $-\varepsilon_{2}=-1+\varepsilon_{6}$ which results from taking 
augmentation, it follows that $\varepsilon_{6}=1$, and $\varepsilon_{2}=0$. 
Thus $u$ is rationally conjugate to a group element, by
\cite[Theorem~2.5]{MaRiSeWe:87} and Proposition~\ref{primeord}.
\end{proof}

The last result is as follows. 
\begin{proposition}\label{pord}
Let $G=\text{\rm PSL}(2,p^{f})$. Then for any $p$-regular torsion unit $u$ 
in $\NU{\ZZ G}$, there is an element of $G$ of the same order as $u$.
\end{proposition}
\begin{proof}
Let $u$ be a $p$-regular torsion unit in $\NU{\ZZ G}$.
If $u$ is of prime power order the claim follows from \cite[(7.3)]{Seh:93}.
Thus we can assume that $u$ has order $rs$ for different primes $r$ and $s$
such that $r\mid p^{f}-1$ and $s\mid p^{f}+1$,
and have to reach a contradiction. 

By Proposition~\ref{primeord} and Remark~\ref{remi2}, we have
\[ \Theta_{3}(u^{r})\sim\text{\rm diag}(\zeta_{s}^{},1,\zeta_{s}^{-1}), \quad
\Theta_{3}(u^{s})\sim\text{\rm diag}(\zeta_{r}^{},1,\zeta_{r}^{-1}) \]
for a primitive $s$-th root of unity $\zeta_{s}$ and a 
primitive $r$-th root of unity $\zeta_{r}$.
This shows that $\Theta_{3}(u)$ has a primitive $rs$-th root of unity
as an eigenvalue, and replacing the roots of unity by suitable powers, 
if necessary, we either have 
$\Theta_{3}(u)\sim\text{\rm diag}
(\zeta_{r}^{}\zeta_{s}^{},\zeta_{r}^{-1},\zeta_{s}^{-1})$ or
$\Theta_{3}(u)\sim\text{\rm diag}
(\zeta_{r}^{}\zeta_{s}^{},1,\zeta_{r}^{-1}\zeta_{s}^{-1})$.
By \cite[(7.3)]{Seh:93}, $G$ has elements $x$ and $y$ of order
$r$ and $s$, respectively. By Theorem~\ref{MT} and Theorem~\ref{T8}, 
these elements can be chosen such that
\begin{align}\label{ag1}
\varphi_{3}(u) & =
1+\sum_{i\in I}\varepsilon_{x^{i}}(u)(\zeta_{r}^{i}+\zeta_{r}^{-i})
+\sum_{j\in J}\varepsilon_{y^{j}}(u)(\zeta_{s}^{j}+\zeta_{s}^{-j})
\end{align}
where the $x^{i}$ ($i\in I$) are representatives of the conjugacy classes
of elements of order $r$, and the $y^{j}$ ($j\in J$) are representatives of 
the conjugacy classes of elements of order $s$. In particular, 
$\varphi_{3}(u)$ is a real number, which shows that only
$\Theta_{3}(u)\sim\text{\rm diag}
(\zeta_{r}^{}\zeta_{s}^{},1,\zeta_{r}^{-1}\zeta_{s}^{-1})$
is possible, and
\begin{align}\label{ag2}
\varphi_{3}(u)=1+\zeta_{r}^{}\zeta_{s}^{}+\zeta_{r}^{-1}\zeta_{s}^{-1}.
\end{align}

Suppose that $r=2$. Then 
$\varphi_{3}(u)=1-\zeta_{s}^{}-\zeta_{s}^{-1}$, and by (\ref{ag1})
\begin{align*}
\varphi_{3}(u) & = 1-2\varepsilon_{x}(u)
+\sum_{j=1}^{(s-1)/2}\varepsilon_{y^{j}}(u)(\zeta_{s}^{j}+\zeta_{s}^{-j}) \\
& = 1+\sum_{j=1}^{(s-1)/2} (2\varepsilon_{x}(u)+\varepsilon_{y^{j}}(u)) 
(\zeta_{s}^{j}+\zeta_{s}^{-j}),
\end{align*}
so $2\varepsilon_{x}(u)+\varepsilon_{y}(u)=-1$, and
$\varepsilon_{y^{j}}(u)=-2\varepsilon_{x}(u)$ for $1<j\leq(s-1)/2$.
It follows that 
\[ 1=\varepsilon_{x}(u)+\sum_{j=1}^{(s-1)/2}\varepsilon_{y^{j}}(u)= 
-1+(2-s)\varepsilon_{x}(u), \]
so $\varepsilon_{x}(u)=2/(2-s)$. By Proposition~\ref{ord6}, we have 
$s>3$. But then $0\neq\varepsilon_{x}(u)<1$, a contradiction.

Thus we can assume that both $r$ and $s$ are odd.
Let $\tau$ be a Galois automorphism satisfying 
$\zeta_{r}^{\tau}=\zeta_{r}^{-1}$ and $\zeta_{s}^{\tau}=\zeta_{s}$.
By (\ref{ag1}), $\varphi_{3}(u)^{\tau}=\varphi_{3}(u)$, so
$\zeta_{r}^{-1}\zeta_{s}^{}+\zeta_{r}^{}\zeta_{s}^{-1}=
\zeta_{r}^{}\zeta_{s}^{}+\zeta_{r}^{-1}\zeta_{s}^{-1}$ by (\ref{ag2}),
which is impossible. The proposition is proved.
\end{proof}

In concluding, we record that the Zassenhaus conjecture is verified 
for the groups $\text{\rm PSL}(2,p)$, $p=7,11,13$, 
by Propositions~\ref{PSLa}, 
\ref{PSLb}, \ref{primeord}, \ref{ord6} and \ref{pord}.



\providecommand{\bysame}{\leavevmode\hbox to3em{\hrulefill}\thinspace}

\end{document}